\documentclass[12pt]{article}
\usepackage{amsthm,amsfonts,amssymb,amscd}

\begin{document}

\begin{center}
\large \bf Birational geometry of singular Fano varieties
\end{center}\vspace{0.5cm}

\centerline{A.V.Pukhlikov}\vspace{0.5cm}

\parshape=1
3cm 10cm \noindent {\small \quad\quad\quad \quad\quad\quad\quad
\quad\quad\quad {\bf }\newline We prove divisorial canonicity of
Fano hypersurfaces and double spaces of general position with
elementary singularities.} \vspace{1.5cm}

\section*{Introduction}\vspace{0.1cm}

{\bf 0.1. The theorem on Fano direct products.} Recall [1] that a
primitive Fano variety $F$ (that is, a smooth Fano variety with
$\mathop{\rm Pic}F={\mathbb Z}K_F$) satisfies the condition of
{\it divisorial canonicity}, or the condition ($C$) (respectively,
the condition of {\it divisorial log canonicity}, or the
condition ($L$)), if for any effective divisor $D\in|-nK_F|$,
$n\geq 1$, the pair
\begin{equation}\label{i1}
(F,\frac{1}{n}D)
\end{equation}
has canonical (respectively, log canonical) singularities. If the
pair (\ref{i1}) has canonical singularities for a general divisor
$D\in \Sigma \subset|-nK_F|$ of any {\it movable} linear system
$\Sigma$, then we say that $F$ satisfies the condition of {\it
movable canonicity}, or the condition ($M$).

Explicitly, the condition ($C$) is formulated in the following
way: for any birational morphism $\varphi\colon \widetilde F\to
F$ and any exceptional divisor $E\subset \widetilde F$ the
following inequality
\begin{equation}\label{i2}
\nu_E(D)\leq na(E)
\end{equation}
holds. The inequality (\ref{i2}) is opposite to the classical {\it
Noether-Fano} inequality, see [2] and the bibliography in that
paper. The condition ($L$) is weaker: the inequality
\begin{equation}\label{i3}
\nu_E(D)\leq n(a(E)+1)
\end{equation}
is required. In (\ref{i2}) and (\ref{i3}) the number $a(E)$ is the
discrepancy of the exceptional divisor $E\subset \widetilde F$
with respect to the model $F$. The inequality (\ref{i3}) is
opposite to the {\it log Noether-Fano inequality}. The condition
($M$) means that (\ref{i2}) holds for a general divisor $D$ of any
movable linear system $\Sigma\subset |-nK_F|$ and any discrete
valuation $\nu_E$.

In [1] the following fact was shown.

{\bf Theorem 1.} {\it Assume that primitive Fano varieties
$F_1,\dots,F_K$, $K\geq 2$, satisfy the conditions ($L$) and
($M$). Then their direct product
$$
V=F_1\times\dots\times F_K
$$
is a birationally superrigid variety. In particular,

{\rm (i)} Every structure of a rationally connected fiber space
on the variety $V$ is given by a projection onto a direct factor.
More precisely, let $\beta\colon V^{\sharp}\to S^{\sharp}$ be a
rationally connected fiber space and $\chi\colon V\dashrightarrow
V^{\sharp}$ a birational map. Then there exists a subset of
indices
$$
I=\{i_1,\dots,i_k\}\subset \{1,\dots,K\}
$$
and a birational map
$$
\alpha\colon F_I=\prod\limits_{i\in I}F_i \dashrightarrow
S^{\sharp},
$$
such that the diagram
$$
\begin{array}{rcccl}
& V &\stackrel{\chi}{\dashrightarrow} & V^{\sharp}&\\
\pi_I &\downarrow & &\downarrow &\beta\\
& F_I & \stackrel{\alpha}{\dashrightarrow}& S^{\sharp}&
\end{array}
$$
commutes, that is, $\beta\circ\chi=\alpha\circ\pi_I$, where
$$
\pi_I\colon\prod\limits^K_{i=1}F_i\to \prod\limits_{i\in I}F_i
$$
is the natural projection onto a direct factor.

{\rm (ii)} Let $V^{\sharp}$ be a variety with ${\mathbb
Q}$-factorial terminal singularities, satisfying the condition
$$
\mathop{\rm dim}\nolimits_{\mathbb Q}(\mathop{\rm
Pic}V^{\sharp}\otimes{\mathbb Q})\leq K,
$$
and $\chi\colon V\dashrightarrow V^{\sharp}$ a birational map.
Then $\chi$ is a (biregular) isomorphism.

{\rm (iii)} The groups of birational and biregular self-maps of
the variety $V$ coincide:
$$
\mathop{\rm Bir}V=\mathop{\rm Aut}V.
$$
In particular, the group $\mathop{\rm Bir}V$ is finite.

{\rm (iv)} The variety $V$ admits no structures of a fibration
into rationally connected varieties of dimension strictly smaller
than $\mathop{\rm min}\{\mathop{\rm dim}F_i\}$. In particular, $V$
admits no structures of a conic bundle or a fibration into
rational surfaces.

{\rm (v)} The variety $V$ is non-rational.}

For the precise definition of birational (super)rigidity, a
discussion of its properties and a list of examples of
birationally (super)rigid varieties, see [2].

There are few doubts that the condition $(C)$ is strictly sharper
than $(L)$: it is easy to give examples of varieties that do not
satisfy the condition $(C)$ but most probably satisfy $(L)$.
However, the existing techniques makes it possible (for typical
Fano varities) to prove $(L)$ via $(C)$ only, showing at once the
stronger property (which automatically imply $(L)$ and $(M)$.)

Note also that the smoothness of the variety $F$ is absolutely
inessential: ${\mathbb Q}$-factoriality is sufficient. If the pair
$(F,\frac{1}{n}D)$ is canonical for any effective divisor
$D\sim-nK_F$ (the equivalence up to multiplication by some
positive integer), then the variety $F$ is covered by Theorem 1.
In [1] this theorem is formulated and proved in the assumption of
smoothness for the only reason that it is then applied to smooth
varieties. The proof given in [1] works word for word for
${\mathbb Q}$-factorial Fano varieties, satisfying the conditions
$(L)$ and $(M)$. This point will be meant in the sequel without
special comments.\vspace{0.3cm}

%%%%%%%%%%%%%%%%%%%%%%%%%%%%%%%%%%%%%%%%%%%%%%%%%%%%%%%%%%%%%%%%%%%
%%%%%%%%%%%%%%%%%%%%%%%%%%%%%%%  subsection 0.2

{\bf 0.2. Examples of divisorially canonical varieties.} In [1]
the divisorial canonicity was shown for two classes of varieties:
generic double spaces of dimension three and higher and generic
Fano hypersurfaces of index 1 and dimension $\geq 5$. Recall what
is meant by being generic.

First let
$$
F\stackrel{\sigma}{\to}{\mathbb P}^M
$$
be a Fano double space branched over a smooth hypersurface
$W=W_{2M}\subset {\mathbb P}^M$ of degree $2M$, $M\geq 3$. For a
point $x\in W$ fix a system of affine coordinates $z_1,\dots,z_M$
on ${\mathbb P}^M$ with the origin at $x$ and set
$$
w= q_1+ q_2+\dots + q_{2M}
$$
to be the equation of the hypersurface $W$, $q_i= q_i(z_*)$ are
homogeneous polynomials of degree $\mathop{\rm deg}q_i=i$. One
should consider the three cases $M\geq 5$, $M=4$ and $M=3$
separately. For convenience of notations assume that $q_1\equiv
z_1$. Set also
$$
\bar q_i=\bar
q_i(z_2,\dots,z_M)=q_i|_{\{z_1=0\}}=q_i(0,z_2,\dots,z_M).
$$
For $M\geq 5$ we say that the Fano double space $F$ is regular at
the point $x$, if the rank of the quadratic form $\bar q_2$ is at
least 2.

Assume that $M=3$ or 4. In that case we require that the quadratic
form $\bar q_2$ is non-zero and moreover

(i) either $\mathop{\rm rk}{\bar q}_2\geq 2$ (as above),

(ii) or $\mathop{\rm rk}{\bar q}_2=1$ and the following
additional condition is satisfied. Without loss of generality we
assume in this case that
$$
\bar q_2= z^2_2.
$$
Now for $M=4$ we require that the following cubic polynomial in
the variable $t$,
$$
\bar q_3(0,1,t)= q_3(0,0,1,t)
$$
has three distinct roots.

For $M=3$ we require that at least one of the following two
polynomials in the variable $z_3$,
$$
\bar q_3(0,z_3)\quad \mbox{or}\quad \bar q_4(0,z_3)
$$
(they are of the form $\alpha z^m_3$, $\alpha\in{\mathbb C}$,
$m=3,4$) is non-zero.

Let
$$
{\cal W}={\mathbb P}(H^0({\mathbb P}^M,{\cal O}_{{\mathbb
P}^M}(2M)))
$$
be the space of hypersurfaces of degree $2M$. Let ${\cal W}_{\rm
reg}\subset {\cal W}$ be the set of branch divisors, satisfying
the regularity condition at every point. The following fact was
shown in [1].

{\bf Theorem 2.} {\it The set ${\cal W}_{\rm reg}$ is non-empty.
For any branch divisor $W\in {\cal W}_{\rm reg}$ the
corresponding Fano double space $\sigma\colon F\to{\mathbb P}^M$,
branched over $W$, satisfies the condition $(C)$.}

Now let
$$
F=F_{M+1}\subset{\mathbb P}={\mathbb P}^{M+1}
$$
be a smooth Fano hypersurface of degree $M+1$, $M\geq 5$. For a
point $x\in F$ fix a system of affine coordinates
$z_1,\dots,z_{M+1}$ with the origin at $x$. Let
$$
f=q_1+q_2+\dots+q_{M+1}
$$
be the equation of the hypersurface $F$, $q_i=q_i(z_*)$ are
homogeneous polynomials of degree $\mathop{\rm deg} q_i=i$. Set
$$
f_i=q_1+\dots+q_i
$$
to be the left segments of the polynomial $f$, $i=1,\dots,M$. Let
us formulate the regularity conditions.

(R1.1) The sequence
$$
q_1,q_2,\dots,q_{M}
$$
is regular in the ring ${\cal O}_{x,{\mathbb P}}$, that is, the
system of equations
$$
q_1=q_2=\dots=q_{M}=0
$$
defines a one-dimensional subset, a finite set of lines in
${\mathbb P}$, passing through the point $x$.

(R1.2) The linear span of any irreducible component of the closed
algebraic set
$$
q_1=q_2=q_3=0
$$
in ${\mathbb C}^{M+1}$ is the hyperplane $q_1=0$ (that is, the
tangent hyperplane $T_xF$).

(R1.3) The closed algebraic set
$$
\overline{\{f_1=f_2=0\}\cap F}=\overline{\{q_1=q_2=0\}\cap
F}\subset {\mathbb P}
$$
(the bar $\bar{}$ means the closure in ${\mathbb P}$) is
irreducible and any section of this set by a hyperplane $P\ni x$
is

\begin{itemize}

\item either also irreducible and reduced,

\item or breaks into two irreducible components
$B_1+B_2$, where $B_i=F\cap S_i$ is the section of $F$ by a plane
$S_i\subset{\mathbb P}$ of codimension 3, and moreover
$\mathop{\rm mult}\nolimits_x B_i=3$,

\item or is non-reduced and is of the form $2B$, where
$B=F\cap S$ is the section of $F$ by a plane $S$ of codimension 3,
and moreover $\mathop{\rm mult}\nolimits_x B=3$.

\end{itemize}

Set
$$
{\cal F}={\mathbb P}(H^0({\mathbb P},{\cal O}_{\mathbb P}(M+1)))
$$
to be the space of hypersurfaces of degree $M+1\geq 6$. Let
${\cal F}_{\rm reg}\subset{\cal F}$ be the set of Fano
hypersurfaces, satisfying the conditions (R1.1-R1.3) at every
point (in particular, every hypersurface $F\in {\cal F}_{\rm
reg}$ is smooth). The next fact was shown in [1].

{\bf Theorem 3.} {\it The set ${\cal F}_{\rm reg}$ is non-empty.
Every Fano hypersurface $F\in{\cal F}_{\rm reg}$ satisfies the
condition $(C)$.} \vspace{0.3cm}

%%%%%%%%%%%%%%%%%%%%%%%%%%%%%%%%%%%%%%%%%%%%%%%%%%%%%%%%%%%%%%%%%
%%%%%%%%%%%%%%%%%%%%%%%%%%%%%%   subsection 0.3

{\bf 0.3. Formulation of the main result.} The aim of this paper
is to show that a Fano variety of any of the two types considered
above still satisfies the condition $(C)$ when it aquires
elementary singularities.

{\bf Theorem 4.} {\it The double space $\sigma\colon V\to{\mathbb
P}^M$, $M\geq 3$, branched over a hypersurface $W\subset{\mathbb
P}^M$ of degree $2M$ with isolated non-degenerate quadratic
singularities, satisfying the regularity conditions of Theorem 2
at every smooth point, satisfies the condition $(C)$.}

{\bf Proof} is not hard: the theorem follows immediately from
Proposition 1.3 (and the fact that the multiplicity of an
irreducible subvariety $Y$ at every poiny does not exceed its
anticanonical degree $(Y\cdot(-K_V)^{\rm dim Y})$. Using the
method of [1], it is easy to check that a generic hypersurface
$W\in{\cal W}$ with a fixed singular point $o\in{\mathbb P}^M$ is
regular at every smooth point, so that varieties described in
Theorem 4 exist.

Now let us consider Fano hypersurfaces
$$
V=V_{M+1}\subset{\mathbb P}={\mathbb P}^{M+1}.
$$
We assume that $M=\mathop{\rm dim}V \geq 8$. This restriction is
connected with the techniques of the proof of Theorem 5: for the
smaller values of $M$ the arguments do not work. Let us formulate
the regularity condition for a non-degenerate double point $o\in
V$.

Let $z_1,\dots,z_{M+1}$ be a system of affine coordinates on the
space ${\mathbb P}$ with the origin at the point $o$,
$$
f=q_2+q_3+\dots+q_{M+1}=0
$$
the equation of the hypersurface $V$, decomposed into homogeneous
components, where $q_2(z_*)$ is a non-degenerate quadratic form.
We say that $V$ is {\it regular} at the point $o$, if the
following conditions hold:

(R2.1) the sequence $q_2,\dots,q_{M+1}$ is regular in ${\cal
O}_{o,{\mathbb P}}$, that is, the system of equations
$$
q_2=\dots=q_{M+1}=0
$$
defines a finite set of points in ${\mathbb P}^M$ (corresponding
to the lines in ${\mathbb P}$, passing through the point $o$ and
lying on $V$),

(R2.2) for any $k\in\{2,3,4,5\}$ and any linear subspace
$P\subset{\mathbb P}$ of codimension two, containing the point
$o$, the closed algebraic set
\begin{equation}\label{i4}
V\cap P\cap\{q_2=0\}\cap\dots\cap\{q_k=0\}
\end{equation}
is irreducible and has multiplicity precisely $(k+1)!$ at the
point $o$.

{\bf Remark 0.1.} The condition that the multiplicity of the set
(\ref{i4}) at the point $o$ is $(k+1)!$, in terms of the
polynomials $q_i$ means that the intersection of the closed set
$$
\{q_2=0\}\cap\dots\cap\{q_{k+1}=0\}\subset{\mathbb P}^M
$$
with any linear subspace of codimension two is of codimension
$k+2$ in ${\mathbb P}^M$ (and of degree $(k+1)!$).

The methods of [1] combined with the methods of [3] make it
possible to show that a generic hypersurface $V\subset{\mathbb
P}^M$ with a fixed double point $o\in{\mathbb P}$ satisfies the
conditions (R1.1-R1.3) at every point $x\neq o$, $x\in V$. The
main fact, that is, that the condition (R1.1) holds, is shown in
[4]. The additional conditions (R1.2-R1.3) are checked directly:
the presence of a fixed singularity does not affect the arguments
of [1, Sec. 2.3]. Finally, the fact that a generic hypersurface
$V\ni o$, singular at the point $o$, satisfies at that point the
conditions (R2.1) and (R2.2), is obvious. Therefore, a generic
hypersurface $V\in{\cal F}$ with a fixed double point
$o\in{\mathbb P}$ is regular at every point, in the sense of the
conditions (R1.1-R1.3) or (R2.1-R2.2).

The main result of the paper is

{\bf Theorem 5.} {\it Assume that the Fano hypersurface
$V\subset{\mathbb P}$ of degree $M+1\geq 9$ is regular at every
point, smooth or regular. Then the variety $V$ satisfies the
condition $(C)$.} \vspace{0.3cm}

%%%%%%%%%%%%%%%%%%%%%%%%%%%%%%%%%%%%%%%%%%%%%%%%%%%%%%%%%%%%%%%%%
%%%%%%%%%%%%%%%%%%%%%%%%%%    subsection 0.4

{\bf 0.4. The structure of the paper and the scheme of the
proof.} Let $X$ be an algebraic variety, $D$ an effective
${\mathbb Q}$-divisor, $S\subset X$ an irreducible subvariety
which is not contained entirely in the set of singular points
$\mathop{\rm Sing}X$. We say that $S$ is an {\it isolated centre
of a non (log) canonical singularity of the pair} $(X,D)$, if
there exists a non (log) canonical singularity of this pair, the
centre of which on $X$ is $S$ (that is, for some resolution
$\varphi\colon X^+\to X$ and a prime exceptional divisor
$E\subset X^+$ the Noether-Fano inequality
$$
\nu_E(D)>a(E)
$$
or, respectively, in the log version,
$$
\nu_E(D)>a(E)+1,
$$
holds, where $\varphi(E)=S$), and there are no non (log)
canonical singularities of that pair, the centre of which on $X$
strictly contain $S$. The main technical tool that we use in this
paper to investigate non (log) canonical singularities is the
following

{\bf Proposition 0.1.} {\it Assume that $\mathop{\rm codim}S\geq
2$ and $S$ is an isolated centre of non-canonical singularity of
the pair $(X,D)$. Let $\mu\colon\widetilde X\to X$ be the blow up
of a point $x\in S$ of general position, $E_x\subset\widetilde X$
the exceptional divisor. For some hyperplane $B\subset E_x$ the
inequality
$$
\mathop{\rm mult}\nolimits_xD+\mathop{\rm
mult}\nolimits_B\widetilde D> 2
$$
holds, where $\widetilde D\subset\widetilde X$ is the strict
transform of the divisor $D$ on $\widetilde X$.}

{\bf Proof} in the case when $S=x$ is a non-singular point on $X$
is given in [1, Proposition 3]. The general case reduces to this
one by restricting the pair $(X,D)$ onto a generic smooth germ
$R\ni x$ of dimension $\mathop{\rm codim}S$. The proof is
complete.

Now let $D\sim -nK_V$ be an effective divisor on a variety $V$ of
any of the two types considered in this paper. Theorems 4 and 5
assert that the pair $(V,\frac{1}{n}D)$ is canonical. Assume that
this is not the case. By Proposition 0.1 and the facts proven in
[1], the centre of a non-canonical singularity of this pair can
be a singular point $o\in V$ only. This case should be excluded.

In \S 1 we prove (in the most general assumptions, not using the
regularity conditions) the inequality
$$
\mathop{\rm mult}\nolimits_oD>2n,
$$
which makes it possible to prove at once Theorem 4. In \S 2 we
carry out a local investigation of the pair
$(V^+,\frac{1}{n}D^+)$, where $V^+$ is the blow up of the point
$o$. The main result of \S 2 is the existence of a hyperplane
section of the exceptional quadric, which is of high multiplicity
with respect to the divisor $D^+$ (the strict transform of $D$ on
$V^+$). This makes it possible \S 3 to prove Theorem 5 by
restricting the pair $(V,\frac{1}{n}D)$ onto the hyperplane
section of the variety $V$, corresponding to the hyperplane
section of the exceptional quadric, which was found in \S 2. One
has to perform the operation of restricting onto a hyperplane
section twice. The principal method of proving Theorem 5 is the
method of hypertangent divisors using the regularity conditions
(R2.1-R2.2), see [2, Chapter III] and the bibliography in that
paper.\vspace{0.3cm}

%%%%%%%%%%%%%%%%%%%%%%%%%%%%%%%%%%%%%%%%%%%%%%%%%%%%%%%%%%%%%%%%%%
%%%%%%%%%%%%%%%%%%%%%%%%%%%% subsection 0.5

{\bf 0.5. Remark.} Up to this day, no examples of singular Fano
varieties, satisfying the condition of divisorial canonicity,
were known in the dimension $\geq 4$. The examples of Theorems 4
and 5 are the first ones. In dimension three some examples were
known [5]: they are weighted Fano hypersurfaces; however, since
their anticanonical degree is small, their study was not
difficult.

%%%%%%%%%%%%%%%%%%%%%%%%%%%%%%%%%%%%%%%%%%%%%%%%%%%%%%%%%%%%%%%%%%%
%%%%%%%%%%%%%%%%%%%%%%%%%%%%%%%%%%%%%%%%%%%%%%%%%%%%%%%%%%%%%%%%%%%
%%%%%%%%%%%%%%%%%%%%%%%%%%%%%%         section 1

\section{The connectedness principle and its first applications}

{\bf 1.1. The connectedness principle.} Let $X$, $Z$ be normal
varieties or analytic spaces and $h\colon X\to Z$ a proper
morphism with connected fibers, and $D=\sum d_iD_i$ a $\mathbb
Q$-divisor on $X$.

{\bf Theorem 6 (connectedness principle, [6 ,Theorem 17.4])} {\it
Assume that $D$ is effective $(d_i\geq 0)$ and the class
$$
-(K_X+D)
$$
is $h$-numerically effective and $h$-big. Let
$$
f\colon Y\stackrel{h}{\to}X\stackrel{h}{\to}Z
$$
be a resolution of singularities of the pair $(X,D)$. Set
$$
K_Y=g^*(K_X+D)+\sum e_iE_i.
$$
The support of the ${\mathbb Q}$-divisor
$\sum\limits_{e_i\leq-1}e_iE_i$, that is, the closed algebraic set
$$
\bigcup\limits_{e_i\leq-1}E_i,
$$
is connected in a neighborhood of any fiber of the morphism $f$.}

{\bf Proof} see in [6, Ch. 17].

As an application of the connectedness principle, consider a germ
$o\in V$ of an isolated terminal singularity with the following
properties. Let
$$
\varphi\colon V^+\to V
$$
be the blow up of the point $o$, $E=\varphi^{-1}(o)$ the
exceptional divisor, which is irreducible and reduced. The
varieties $V$, $V^+$ and $E$ have ${\mathbb Q}$-factorial
terminal singularities. Let $\delta=a(E,V)$ be the discrepancy of
$E$, $D$ an effective ${\mathbb Q}$-divisor on $V$, $D^+$ its
strict transform on $V^+$. Define the number $\nu_E(D)$ by the
formula
$$
\varphi^*D=D^++\nu_E(D)E.
$$

{\bf Proposition 1.1.} {\it Assume that the pair $(V,D)$ is not
canonical at the point $o$, which is an isolated centre of a
non-canonical singularity of this pair. Assume also that for some
integer $k\geq 1$ the inequality
\begin{equation}\label{a1}
\nu_E(D)+k\leq\delta
\end{equation}
holds. Then the pair $(V^+,D^+)$ is not log canonical and there
is a non log canonical singularity $\widetilde E\subset\widetilde
V$ of that pair (where $\widetilde V\to V^+$ is some model), the
centre of which
$$
\mathop{\rm centre}(\widetilde E,V^+)\subset E
$$
is of dimension $\geq k$.}

{\bf Proof.} Assuming $V\subset{\mathbb P}^N$ to be projectively
embedded, consider a {\it generic} linear subspace
$P\subset{\mathbb P}^N$ of codimension $k$, containing the point
$o$. Let $\Lambda_P$ be the linear system of hyperplanes,
containing $P$, and $\Lambda$ be the corresponding linear system
of sections of the variety $V$. Let $\varepsilon>0$ be a
sufficiently small rational number of the form $\frac{1}{K}$ and
$$
\{ H_I\,|\,i\in I\}\subset\Lambda
$$
a set of $\sharp I=Kk$ generic divisors. Set
$$
R=D+\sum_{i\in I}\varepsilon H_i,
$$
and let $R^+$ be the strict transform of $R$ on $V^+$.

Obviously, the pair $(V^+,D^+)$ is not log canonical. The centre
of any of its non log canonical singularities is contained in
$E$. Furthermore, being non log canonical is an open property, so
that, slightly decreasing the coeffients in $D$, we may assume
that the strict version of the inequality (\ref{a1}) holds, that
is, $\nu_E(D)+k<\delta$ (whereas other assumptions still hold).

Now consider the pair $(V^+,R^+)$ (we still assume that $V\ni o$
is a germ, so that all constructions are local in a neighborhood
of the point $o$). It is non log canonical, and all its non log
canonical singularities are non log canonical singularities of
the pair $(V^+,D^+)$, with the exception of one additional
singularity, the germ $(P\cap V)^+$ of the section of $V$ by the
plane $P$, that is, the base set of the system $\Lambda$. By the
strict version of the inequality (\ref{a1}), the class
$-(K_{V^+}+R^+)$ is obviously $\varphi$-nef and $\varphi$-big, so
that, applying the connectedness principle (to $X=V^+$, $Z=V$,
$h=\varphi$, $D=R^+$), we conclude that the union of the centres
of non log canonical singularities of the pair $(V^+,R^+)$ on
$V^+$ is connected. Since $P$ is generic, this is only possible if
$(P\cap V)^+$ intersects some centre of a non log canonical
singularity of the pair $(V^+,D^+)$, which should be of dimension
at least $k$. Q.E.D. for the proposition.

The fact which we have just proven will be applied to our case of
a hypersurface singularity $o\in V$ with a smooth exceptional
divisor. \vspace{0.3cm}

%%%%%%%%%%%%%%%%%%%%%%%%%%%%%%%%%%%%%%%%%%%%%%%%%%%%%%%%%%%%%%%%%%%%%%
%%%%%%%%%%%%%%%%%%%%%%%%%%%%%%%%%%  subsection 1.2

{\bf 1.2. Singularities of pairs on a smooth hypersurface.} Let
$X\subset{\mathbb P}^N$ be a smooth hypersurface of degree
$m\in\{2,\dots,N-1\}$, $D\in\,|\,lH_X\,|$ an effective divisor,
which is cut out on $X$ by a hypersurface of degree $l\geq 1$.
(So that $H_X$ is the class of a hyperplane section of $X$.) The
following fact and its proof are well known [3,7].

{\bf Proposition 1.2.} {\it For any $n\geq l$ the pair
$$
(X,\frac{1}{n}D)
$$
is log canonical.}.

{\bf Proof.} We may consider the case $n=l$. Assume the converse:
the pair $(X,\frac{1}{n}D)$ is not log canonical. Since for any
curve $C\subset X$ the inequality
$$
\mathop{\rm mult}\nolimits_CD\leq n
$$
holds (see [2,4]), the centre of a non log canonical singularity
of the pair $(X,\frac{1}{n}D)$ can only be a point. Let $x\in X$
be such a point. Consider a general projection $\pi\colon{\mathbb
P}^N\dashrightarrow{\mathbb P}^{N-1}$. Its restriction onto $X$
is a finite morphism $\pi_X\colon X\to{\mathbb P}^{N-1}$ of degree
$m$, which is an analytic isomorphism at the point $x$, and one
may assume that
$$
\pi^{-1}_X(\pi_X(x))\cap\mathop{\rm Supp}D=\{x\}.
$$
This implies that the germ of the pair $(X,\frac{1}{n}D)$ at the
point $x$ and the germ of the pair $({\mathbb
P}^{N-1},\frac{1}{n}\pi(D))$ at the point $\pi(x)$ are
analytically isomorphic. In particular, the point $\pi(x)$ is an
isolated centre of a non log canonical singularity of the pair
$({\mathbb P}^{N-1},\frac{1}{n}\pi(D))$. However, this is
impossible.

Being non log canonical is an open property, so that for a
rational number $s<n^{-1}$, sufficiently close to $n^{-1}$, the
pair
$$
({\mathbb P}^{N-1},s\pi(D))
$$
still has the point $\pi(x)$ as an isolated centre of a non log
canonical singularity. Let $P\subset{\mathbb P}^{N-1}$ be a
hyperplane, not containing the point $\pi(x)$. By the inequality
$$
smn+1<N
$$
the ${\mathbb Q}$-divisor $-(K_{{\mathbb P}^{N-1}}+s\pi(D)+P)$ is
ample, so that one may apply to the pair
$$
({\mathbb P}^{N-1},s\pi(D)+P)
$$
the connectedness principle of Shokurov and Koll\' ar (in the
notations of Theorem 6, $X={\mathbb P}^{N-1}$, $Z$ is a point,
for the ${\mathbb Q}$-divisor $D$ we take $s\pi(D)+P$, the
conditions of Theorem 6 are satisfied in a trivial way by what
was said above) and obtain a contradiction: the point $\pi(x)$ is
an isolated centre of a non log canonical singularity and the
divisor $P$ comes into the $\mathbb Q$-divisor $s\pi(D)+P$ with
the coefficient one, however $\pi(x)\not\in P$, so that the
connectedness is violated. Q.E.D. for Proposition 1.2.

The fact which we have obtained makes it possible to start a
serious discussion of an isolated terminal singularity $o\in V$
(the naive arguments, even with strong conditions of general
position for the germ $o\in V$, give too weak estimates for the
multiplicity at the point $o$), but not more than that.
\vspace{0.3cm}

%%%%%%%%%%%%%%%%%%%%%%%%%%%%%%%%%%%%%%%%%%%%%%%%%%%%%%%%%%%%%%
%%%%%%%%%%%%%%%%%%%%%%%%%%%%%%%    subsection 1.3

{\bf 1.3. The weak local inequality.} Let $o\in V$ be a germ of
isolated hypersurface terminal singularity. More precisely, if
$\varphi\colon V^+\to V$ is the blow up of the point $o$,
$\varphi^{-1}(o)=E\subset V^+$ is the exceptional divisor, we
assume that $V^+$ and $E$ are smooth, whereas $E$ is isomorphic
to a smooth hypersurface of degree $\mu=\mathop{\rm mult}_oV$ in
${\mathbb P}^M$.

Furthermore, let $D\ni o$ be a germ of a prime divisor,
$D^+\subset V^+$ its strict transform, $D^+\sim-\nu E$ for
$\nu\in{\mathbb Z}_+$, so that the equality
$$
\mathop{\rm mult}\nolimits_oD=\mu\nu
$$
holds.

{\bf Proposition 1.3.} {\it Assume that the pair
$(V,\frac{1}{n}D)$ is not canonical at the point $o$, which is an
isolated centre of a non-canonical singularity of that pair. Then
the inequality}
\begin{equation}\label{a2}
\nu>n
\end{equation}
{\it holds.}

{\bf Proof.} Assume the converse: $\nu\leq n$. Then the pair
$(V^+,\frac{1}{n}D^+)$ is not canonical, and moreover, the centre
of any non-canonical singularity of this pair (that is, of any
maximal singularity of the divisor $D^+$) is contained in the
exceptional divisor $E$. By the inversion of adjunction the pair
$(E,\frac{1}{n}D^+_E)$, where $D^+_E=D^+\,|\,_E$, is not log
canonical. Let $H_E=-E\,|\,_E$ be the generator of the Picard
group $\mathop{\rm Pic}E$, that is, the hyperplane section of $E$
with respect to the embedding $E\subset{\mathbb P}^M$. We get
$$
D^+_E\sim-\nu E\,|\,_E=\nu H_E.
$$
Since $\nu\leq n$, the non log canonicity of the pair
$(E,\frac{1}{n}D^+_E)$ contradicts to Proposition 1.2. Q.E.D.

{\bf Remark 1.1.} Not using inversion of adjunction, the best one
can get using explicit geometric methods, even with the
conditions of general position for $E$, is the inequality
$\nu>\frac{n}{2}$, which is much weaker than (\ref{a2}). But the
inequality (\ref{a2}) is also insufficient for excluding maximal
singularities on the typical Fano varieties.

{\bf Proof of Theorem 4.} Assume that the pair $(V,\frac{1}{n}D)$
is not canonical. By what was proven in [1], only a singular point
$o$ can be the centre of a non-canonical singularity. The divisor
$D\sim-nK_V$ can be assumed to be irreducible. By Proposition
1.3, we get the inequality
$$
\mathop{\rm mult}\nolimits_oD>2n,
$$
however the anticanonical degree of the divisor $D$ is
$$
\mathop{\rm deg}D=(D\cdot(-K_V)^{\rm dim V-1})=2n.
$$
Therefore, $\mathop{\rm mult}_oD>\mathop{\rm deg}D$, which is
impossible (see a detailed discussion in [2, Chapter II]). This
contradiction proves the theorem.

%%%%%%%%%%%%%%%%%%%%%%%%%%%%%%%%%%%%%%%%%%%%%%%%%%%%%%%%%%%%%%%%%%%%%
%%%%%%%%%%%%%%%%%%%%%%%%%%%%%%%%%%%%%%%%%%%%%%%%%%%%%%%%%%%%%%%%%%%%%
%%%%%%%%%%%%%%%%%%%%%%%%%%%%%%% section 2

\section{A local investigation of a divisor at a quadratic point}

{\bf 2.1. Effective divisors on quadrics.} Let $Q\subset{\mathbb
P}^M$ be a non-degenerate quadric, $H_Q\in\mathop{\rm Pic}Q$ the
class of a hyperplane section, $B\subset Q$ an irreducible
subvariety.

{\bf Definition 2.1.}  We say that an effective divisor $D$ on
$Q$ satisfies {\it the condition} $H(n)$ ({\it with respect to}
$B$), where $n\geq 1$ is a fixed integer, if for any point of
general position $p\in B$ there exists a hyperplane $F(p)\subset
E_p$ in the exceptional divisor $E_p=\varphi^{-1}_p(p)$ of the
blow up $\varphi_p\colon Q_p\to Q$ of the point $p$, for which
the inequality
\begin{equation}\label{b1}
\mathop{\rm mult}\nolimits_pD+\mathop{\rm
mult}\nolimits_{F(p)}\widetilde D> 2n,
\end{equation}
holds, where $\widetilde D\subset Q_p$ is the strict transform of
the divisor $D$.

Note that the divisor $D$ is not assumed to be irreducible and
the number $n$ does not depend on the point $p$. The hyperplane
$F(p)$ depends algebraically on the point $p$; this is assumed
everywhere in the sequel without special comments. Let $l\geq 1$
be the degree of a hypersurface in ${\mathbb P}^M$ that cuts out
$D$ on $Q$, that is,
$$
D\sim lH_Q.
$$

{\bf Proposition 2.1.} {\it  Assume that $\mathop{\rm dim}B\geq
3$. Assume furthermore that the effective divisor $D$ satisfies
the condition $H(n)$ with respect to $B$. Then the following
alternative takes place: either

1) the inequality $l> 2n$ holds (this case will be referred to as
the simple one), or

2) there exists a hyperplane section $Z\subset Q$, containing
entirely the subvariety $B$, such that for a point of general
position $p\in B$ in the notations above
$$
F(p)=\widetilde Z\cap E_p,
$$
where $\widetilde Z\subset Q_p$ is the strict transform, whereas
$Z$ comes into the divisor $D$ with the multiplicity
$$
a>2n-l
$$
(that is, the presentation $D=aZ+D^*$ takes place, where $D^*$
does not contain $Z$ as a component; this case we will call the
hard one).}

{\bf Proof.} Let us assume at once that the simple case does not
realize, that is, $l\leq 2n$. In the notations above for a point
of general position $p\in B$ let $\Lambda_p\subset |H_Q|$ be the
pencil of hyperplane sections, the strict transform of which
$\widetilde\Lambda_p$ on $Q_p$ cuts out on $E_p$ the hyperplane
$F(p)$, that is,
$$
\widetilde\Lambda_p\cap E_p=F(p).
$$
The exceptional divisor $E_p$ is the projectivization of the
tangent space $T_pQ\cong{\mathbb C}^{M-1}$. Let
$$
[F(p)]\subset T_pQ
$$
be the hyperplane, the projectivization of which is $F(p)$.
Consider now $T_pQ$ as an {\it embedded} tangent space (in some
affine chart ${\mathbb C}^M\subset{\mathbb P}^M$) and let
$\overline{T_pQ}\subset{\mathbb P}^M$ be its closure (that is,
the hyperplane in ${\mathbb P}^M$, tangent to $Q$ at the point
$p$). Respectively, let
$$
\overline{F(p)}=\overline{[F(p)]}\subset\overline{T_pQ}
$$
be the closure of the subspace $[F(p)]$. This is a linear
subspace in ${\mathbb P}^M$ of codimension two. It is easy to see
that the base set (and subscheme) of the pencil $\Lambda_p$ is
$$
\mathop{\rm Bs}\Lambda_p =\overline{F(p)}\cap Q.
$$
denote this subset by the symbol $\Theta(p)$. Set
$$
Z=\overline{\bigcup_{p\in B}\Theta(p)}
$$
(where the union is taken over the points of some open subset of
$B$, whereas the overline means the closure).

Note that $\Theta(p)$ is a quadric in $\overline{F(p)}$ with at
least one singular point $p$ and at most a line of double points
(containing $p$). Since $\mathop{\rm dim}B\geq 3$, for a pair of
distinct points of general position $p_1\neq p_2$ we get
$\Theta(p_1)\neq\Theta(p_2)$, which implies that either $Z=Q$ or
$Z$ is a prime divisor on $Q$.

Let $R\in\Lambda_p$ be a general element of the pencil,
$D_R=D\,|\,_R$ the restriction of the divisor $D$. By the
inequality (\ref{b1}) we get
$$
\mathop{\rm mult}\nolimits_pD_R>2n.
$$
The variety $R$ is smooth at the point $p$. We get the
presentation
$$
D_R=a\Theta(p)+D^*_R,
$$
where $D^*_R$ does not contain $\Theta(p)$ as a component,
$a\in{\mathbb Z}_+$ is some non-negative integer. Since
$\Lambda_p$ cuts out $\Theta(p)$ with the multiplicity one, we get
$$
a=\mathop{\rm mult}\nolimits_{\Theta(p)}D.
$$
In particular, if $a\geq 1$, then $Z$ is a divisor on $Q$.

{\bf Lemma 2.1.} {\it The inequality $a\geq 1$ holds.}

{\bf Proof.} Assume the converse: $a=0$. It is easy to see that
$$
\Theta(p)=R\cap T_pR.
$$
The intersection $\Theta(p)\cap D_R$ is of codimension two on
$R$, so that the effective cycle $(\Theta(p)\circ D_R)_R$ is well
defined. Therefore,
$$
2l=\mathop{\rm deg}(\Theta(p)\circ D_R)\geq\mathop{\rm
mult}\nolimits_p(\Theta(p)\circ D_R)>4n,
$$
so that, contrary to the assumption above, the inequality $l>2n$
holds. The contradiction proves the lemma. Q.E.D.

Now we get
$$
D=aZ+D^*,
$$
where the divisor $D^*$ does not contain $Z$ as a component. Set
$$
Z\sim l_ZH_Q,\quad D^*\sim l^*H_Q,
$$
so that the equality
$$
l=al_Z+l^*.
$$
holds.

{\bf Lemma 2.2.} {\it $Z\subset Q$ is a hyperplane section:
$l_Z=1$.}

{\bf Proof.} Since $Z$ is a prime divisor, the set
$$
\Delta=\bigcup_{p\in B}\overline{F(p)}
$$
(where the union is taken over an open subset of $B$) can not be
dense in ${\mathbb P}^M$. Since $\mathop{\rm dim}B\geq 3$, in
this union an at least two-dimensional family of linear subspaces
of codimension two is present. By the following elementary lemma,
the closure of $\Delta$ in ${\mathbb P}^M$ is a hyperplane.
Q.E.D. for the lemma.

{\bf Lemma 2.3.} {\it If the surface $S\subset{\mathbb P}^3$
contains a two-dimensional family of lines, then $S$ is a plane.}

{\bf Proof.} There is a one-dimensional family of lines through a
generic point of the surface $S$ (actually, just one point is
sufficient). The case of a cone is obvious. Q.E.D. for the lemma.

By Lemma 2.2, for a general point $p\in B$ we get
$$
\Theta(p)=Z\cap R
$$
(the left-hand side is contained in the right-hand one, however,
in the right-hand side we have a section of $Q$ by a linear
subspace of codimension two). Therefore, $D^*_R=D^*\,|\,_R$.
Obviously,
$$
\mathop{\rm deg}D^+_R=2(l-a)=2l^*.
$$
Arguing as in the proof of Lemma 2.1, that is, intersecting
$D^*_R$ with $Z$ (or with $\Theta(p)$ inside $R$), we obtain the
estimate
$$
2\mathop{\rm mult}\nolimits_pD^*_R\leq\mathop{\rm deg}D^*_R,
$$
which in terms of our integral parameters gives the estimate
$$
2a+l^*=a+l>2n.
$$
The divisor $Z$ contains $B$ and has at most one singular point.
Therefore, for a general point $p\in B$ the divisor $Z$ is
non-singular at $p$ and thus
$$
F(p)=\widetilde Z\cap E_p,
$$
as we claimed. Q.E.D. for Proposition 2.1.

{\bf Remark 2.1.} If $\mathop{\rm dim}B\geq 4$, then the same
arguments, word for word, give the claim of Proposition 2.1 for
the case when $Q$ is a cone over a non-degenerate quadric.
\vspace{0.3cm}

%%%%%%%%%%%%%%%%%%%%%%%%%%%%%%%%%%%%%%%%%%%%%%%%%%%%%%%%%%%%%%%%%%%%
%%%%%%%%%%%%%%%%%%%%%%%%%%%%%%%     subsection 2.2

{\bf 2.2. Hyperplane section of high multiplicity.} Let $o\in V$
be a germ of a non-degenerate quadratic singularity,
$\varphi\colon V^+\to V$ the blow up of the point $o$,
$E=\varphi^{-1}(o)\subset V^+$ the exceptional divisor. Let $D\ni
o$ be a germ of an effective divisor, where the pair
$(V,\frac{1}{n}D)$  has the point $o$ as an isolated centre of a
non-canonical singularity. Let $D^+\subset V^+$ be the strict
transform of $D$ on $V^+$, $D^+\sim-lE$. Assume that $l\leq 2n$,
so that the pair $(V^+,\frac{1}{n}D^+)$ is not log canonical. Let
$S\subset E$ be a centre of a non log canonical singularity,
which has the maximal dimension (in particular, $S$ is an
isolated centre of a non log canonical singularity), $\mathop{\rm
dim} S\geq 3$. In particular, the inequality
\begin{equation}\label{b2}
\mathop{\rm mult}\nolimits_SD^+>n.
\end{equation}
holds. The following fact is true.

{\bf Proposition 2.2.} {\it One of the following two cases takes
place:

{\rm 1)} $S$ is a hyperplane section of the quadric $E$,

{\rm 2)} there is a hyperplane section $Z\supset S$ of the quadric
$E$, satisfying the inequality}
\begin{equation}\label{b3}
\mathop{\rm mult}\nolimits_ZD^+>\frac{2n-l}{3}.
\end{equation}

{\bf Proof.} If $S\subset E$ is a prime divisor, then by the
inequality (\ref{b2}) we get $l>nl_S$, where $S\sim l_SH_E$ and
$H_E=-E\,|\,_E$ is the hyperplane section of the quadric $E$.
Since by assumption $l\leq 2n$, this implies that $l_S=1$, that
is, $S$ is a hyperplane section (case 1)).

Assume that $\mathop{\rm codim}_ES\geq2$, that is, the case 1)
does not realize. Since the pair $(V^+,\frac{1}{n}D)^+$ is not
log canonical at $S$, for a generic point $p\in S$ there is a
hyperplane $\prod(p)\subset E^+_p$ in the exceptional divisor
$E^+_p=\varphi^{-1}_p(p)$ of the blow up
$$
\varphi_p\colon V^+_p\to V^+
$$
of the point $p$, satisfying the inequality
\begin{equation}\label{b4}
\mathop{\rm mult}\nolimits_pD^+\mathop{\rm
mult}\nolimits_{\prod(p)}D^+_p>2n,
\end{equation}
where $D^+_p\subset V^+_p$ is the strict transform with respect
to the blow up of the point $p$. Let $E(p)\subset V^+_p$ be the
strict transform of the exceptional quadric.

{\bf Lemma 2.4.} {\it For a point $p$ of general position}
$$
\prod (p)\neq E(p)\cap E^+_p.
$$

{\bf Proof.} Otherwise for the restriction $D^+_E=(D^+\circ E)$
we get the inequality
$$
\mathop{\rm mult}\nolimits_pD^+_E>2n,
$$
which is true for almost all points $p\in S$. Since the point $p$
runs through a set of positive dimension, for a general point the
divisor $D^+_E$ does not contain the hyperplane section $E\cap
T_pE$ as a component. This, as in the proof of Lemma 2.1,
immediately implies that $l>2n$. This is a contradiction. Q.E.D.
for the lemma.

For a point of general position $p\in S$ set
$$
F(p)=\prod(p)\cap E(p).
$$
The following inequalities hold:
\begin{equation}\label{b5}
\mathop{\rm mult}\nolimits_pD^+\mathop{\rm
mult}\nolimits_{F(p)}D^+_p+>2n,
\end{equation}
which is weaker than (\ref{b4}), and
\begin{equation}\label{b6}
\mathop{\rm mult}\nolimits_pD^+_E+\mathop{\rm
mult}\nolimits_{F(p)}\widetilde D^+_E>2n,
\end{equation}
which follows from (\ref{b5}), where $\widetilde D^+_E\subset
E(p)$ is the strict transform. All geometric objects,
participating in the inequality (\ref{b6}), are defined in terms
of the quadric $E$, that is, they do not require addressing $V^+$
and $V^+_p$.

Let us apply to the divisor $D^+_E$ on the quadric $E$
Proposition 2.1. By our assumptions, the simple case does not
realize, so that there exist a hyperplane section $Z\supset S$,
cutting out the hyperplane $F(p)$ on $E_p=E(p)\cap E^+_p$, that
is, $F(p)=\widetilde Z\cap E_p$. The section $Z$ comes into the
divisor $D^+_E$ with a multiplicity strictly higher than $2n-l$.
In particular, $\mathop{\rm mult}_ZD^+>0$. \vspace{0.3cm}

%%%%%%%%%%%%%%%%%%%%%%%%%%%%%%%%%%%%%%%%%%%%%%%%%%%%%%%%%%%%%%%%
%%%%%%%%%%%%%%%%%%%%%%%%%%%%%%%%% subsection 2.3

{\bf 2.3. Estimating the multiplicity of the hyperplane section.}
Let us prove the inequality (\ref{b3}). Both the assumptions and
all claims of Proposition 2.2 are local at the point $o$. A
non-degenerate quadratic singularity is analytically equivalent
to the germ of the quadric
$$
\{z^2_1+\dots+z^2_{M+1}=0\}\subset{\mathbb C}^{M+1},
$$
so that one can assume that the divisor $D$ is given by an
equation
$$
f=q_l(z_*)+q_{l+1}(z_*)+\dots,
$$
where the quadric $q(z_*)=\sum\limits^{M+1}\limits_{i=1}z^2_i$
divides none of the polynomials $q_i$ (if $q_i\neq 0$). The
affine coordinates $z_*$ can be looked at as homogeneous
coordinates on the exceptional divisor of the blow up of the
origin
$$
o=(0,\dots,0)\in{\mathbb C}^{M+1} ,
$$
and then $\{q=0\}\subset{\mathbb P}^M$ is precisely the
exceptional quadric $E$. The divisor $D^+_E$ is given by the
equation
$$
q_l\,|\,_E=0.
$$
Let $\lambda(z)=0$ be the equation of the hyperplane section
$Z\subset E$. In terms of the coordinates $z_*$ Proposition 2.1
asserts that for some $a_l>2n-l$
$$
q_l=\lambda^{a_l}g+qw,
$$
where $g(z_*)$ and $w(z_*)$ are homogeneous polynomials of the
corresponding degrees. Replacing $q_l$ by $q_l-qw$, we may assume
that $\lambda^{a_l}$ divides the polynomial $q_l$. This implies
that the strict transform $D^+_l$ of the divisor
$$
D_l=\{q_l\,|\,_V=0\}
$$
is of multiplicity $\geq a_l>2n-l$ along $Z$. However, the divisor
$D_l$ is the intersection of the cone $V$ with the cone
$\{q_l=0\}$ (the polynomial $q_l$ is homogeneous), so that the
inequality (\ref{b6}) for $D^+_E=(D^+_l)\,|\,_E$ implies the
inequality (\ref{b5}) for $D^+_l$, where $D^+_l\subset V^+$ is
the strict transform of $D_l$ on $V^+$. Therefore, the following
claim is true:\vspace{0.3cm}

{\it both divisors $D^+$ and $D^+_l$ satisfy the inequality
(\ref{b5}) for a point of general position $p\in
S$.}\vspace{0.3cm}

By linearity of the inequality (\ref{b5}) (and the obvious fact
that
$$
\mathop{\rm mult}\nolimits_pD^+_l=\mathop{\rm
mult}\nolimits_p(D^+_l)\,|\,_E=\mathop{\rm
mult}\nolimits_pD^+_E\geq\mathop{\rm mult}\nolimits_pD^+,
$$
and similarly for $F(p)$), the divisor
\begin{equation}\label{b7}
\varphi^*(f_{l+1}\,|\,_V=0)-lE,
\end{equation}
where
$$
f_{l+1}=q_{l+1}(z_*)+q_{l+2}(z_*)+\dots=f-q_l(z_*),
$$
again satisfies that inequality. Let $k\geq 1$ be the first
index, for which $q_{l+k}(z_*)\,|\,_E\not\equiv 0$,
$$
D_{\geq l+k}=\{f_{l+k}\,|\,_V=0\},
$$
$f_{l+k}=\sum\limits^{\infty}_{i=l+k}q_i(z_*)$. For the strict
transform $D^+_{\geq l+k}$ we get $D^+_{\geq l+k}\sim-(l+k)E$.
Besides, as we mentioned above, the divisor $D^+_{\geq l+k}+kE$
satisfies the inequality (\ref{b5}) at the point of general
position $p\in S$. Since
$$
\mathop{\rm mult}\nolimits_pE=\mathop{\rm
mult}\nolimits_{F(p)}E(p)=1,
$$
the divisor $D^+_{\geq l+k}$ by itself satisfies the inequality
$$
\mathop{\rm mult}\nolimits_pD^+_{\geq l+k}+\mathop{\rm
mult}\nolimits_{F(p)}(D^+_{\geq l+k})_p>2(n-k)
$$
(the lower index $p$ for a divisor means, as usual, the strict
transform on the variety $V^+_p$).

This makes it possible to prove the inequality (\ref{b3}) by
decreasing induction on $l\leq 2n$. The base of induction is the
case $l=2n$: in that case $Z$ comes into $D^+_E$ with a positive
multiplicity, that is, $\mathop{\rm mult}_ZD^+>0$, which is what
we claim in 2) for $l=2n$.

If in the notations above the inequality
$$
l+k\leq 2(n-k)
$$
holds, then by the induction hypothesis we get
$$
\mathop{\rm mult}\nolimits_ZD^+_{\geq
l+k}>\frac{2(n-k)-(l+k)}{3}=\frac{2n-l}{3}-k,
$$
so that the divisor (\ref{b7}), obtained by subtraction from the
equation $f$ the equation of the divisor $D^+_E$, contains $Z$
with a multiplicity strictly higher than
$$
k+\mathop{\rm mult}\nolimits_ZD^+_{\geq l+k}>\frac{2n-l}{3}
$$
(since the divisor (\ref{b7}) contains with the multiplicity $k$
the exceptional quadric $E$), whence, taking into account the
inequality $\mathop{\rm mult}_ZD^+_l>2n-l$, we obtain the
required inequality (\ref{b3}).

If the inequality
$$
l+k>2(n-k),
$$
holds, then we can not apply the induction hypothesis, but in
that case the estimate
$$
k>\frac{2n-l}{3},
$$
holds, so that, arguing as above, we obtain the inequality
(\ref{b3}) all the same, simply because the divisor (\ref{b7})
contains the exceptional quadric $E$ with the multiplicity $k$.
Q.E.D. for Proposition 2.2.

{\bf Remark 2.2.} The claim of Proposition 2.2 remains true for a
germ of an elementary degenerate singularity, when $E$ is a cone
over a non-degenerate quadric, if we assume that $\mathop{\rm
dim}S\geq 4$. No changes in the proof are needed. One can
consider a general hyperplane section of the germ $o\in V$
(containing the point $o$) and apply Proposition 2.2 to that
section.

%%%%%%%%%%%%%%%%%%%%%%%%%%%%%%%%%%%%%%%%%%%%%%%%%%%%%%%%%%%%%%%%%%%
%%%%%%%%%%%%%%%%%%%%%%%%%%%%%%%%%%%%%%%%%%%%%%%%%%%%%%%%%%%%%%%%%%%
%%%%%%%%%%%%%%%%%%%%%%%%%%%%%%% section 3

\section{Exclusion of a maximal singularity}

{\bf 3.1. Estimating the multiplicity at the double point.} Let
$V\subset{\mathbb P}={\mathbb P}^{M+1}$ be a hypersurface of
degree $M+1$, where $M\geq 8$, $o\in V$ an isolated quadratic
singularity, satisfying the regularity conditions (R2.1-R2.2). Set
$$
\varphi\colon V^+\to V
$$
to be the blow up of the point $o$, $E=\varphi^{-1}(o)$ the
exceptional quadric. By the symbol $H$ we denote the class of a
hyperplane section of $V$. We use the notations of Sec. 0.3, in
particular, $q_2=0$ is the equation of the tangent cone at the
point $o$.

Let $D\sim nH$ be an effective divisor, $D^+$ its strict
transform, $D^+\sim nH-\nu E$ for some $\nu\geq 1$.

{\bf Proposition 3.1.} {\it The inequality $\nu\leq\frac32n$
holds.}

{\bf Proof.} This inequality is linear in the divisor $D$, so
that without loss of generality we may assume that $D$ is a prime
divisor. Assume the converse: $\nu>\frac32n$. For the first
hypertangent divisor $D_2=\{q_2\,|\,_V=0\}$ we get
$$
D^+_2\sim 2H-3E,
$$
so that $D_2$ and $D$ are distinct prime divisors. Therefore, the
set-theoretic intersection $D\cap D_2$ is of codimension two and
the effective cycle $Y=(D\circ D_2)$ is well defined. It
satisfies the inequality
$$
\frac{\mathop{\rm mult}_o}{\mathop{\rm
deg}}Y\geq\frac32\cdot\frac{2\nu}{n(M+1)}>\frac{9}{2(M+1)}.
$$
Now let us consider the standard hypertangent systems
$$
\Lambda_k=\left |\sum^{k-2}_{i=o}s_if_{k-i}=0\right |,
$$
where $f_j=q_2+\dots+q_j$ is a left segment of the equation of
$V$ at the point $o$, $s_i$ runs through the set of all
homogeneous polynomials in the coordinates $z_*$ of degree $i$.
By the regularity condition (R2.1), for $r=2,\dots,M$ we get
$$
\mathop{\rm codim}\mathop{\rm Bs}\Lambda_k=k-1,
$$
so that in the usual way [2,4] we construct a sequence of
irreducible subvarieties
$$
Y_2=Y,Y_3,\dots,Y_{M-1}
$$
of codimension $\mathop{\rm codim}Y_i=i$, where $Y_{i+1}\subset
Y_i$ is an irreducible component of the effective cycle
$(Y_i\circ D_{i+2})$, $D_j\in\Lambda_j$ is a hypertangent divisor
of general position, and for $Y_{i+1}$ we take a component with
the maximal ratio $\mathop{\rm mult}_o/\mathop{\rm deg}$. The
effective cycle $(Y_i\circ D_{i+2})$ is well defined, because
$\mathop{\rm codim}\mathop{\rm Bs}\Lambda_{i+2}=i+1$, so that a
general divisor $D_{i+2}$ does not contain $Y_i$. We obtain the
estimate
$$
\frac{\mathop{\rm mult}_o}{\mathop{\rm
deg}}Y_{i+1}\geq\frac{i+3}{i+2}\cdot\frac{\mathop{\rm
mult}_o}{\mathop{\rm deg}}Y_i,
$$
so that for the last subvariety, that is, the curve $Y_{M-1}$, we
obtain the inequality
$$
\frac{\mathop{\rm mult}_o}{\mathop{\rm
deg}}Y_{M-1}\geq\frac{\mathop{\rm mult}_o}{\mathop{\rm
deg}}Y\cdot\frac54\cdot\dots\cdot\frac{M+1}{M}>\frac98,
$$
which is, of course, impossible. The contradiction proves the
proposition. \vspace{0.3cm}

%%%%%%%%%%%%%%%%%%%%%%%%%%%%%%%%%%%%%%%%%%%%%%%%%%%%%%%%%%%%
%%%%%%%%%%%%%%%%%%%%%%%%%     subsection 3.2

{\bf 3.2. Reduction to a hyperplane section.} Let us go back to
the {\bf proof of Theorem 5.} Assume that the pair
$(V,\frac{1}{n}D)$ has the point $o$ as an isolated centre of a
non-canonical singularity. By linearity of the Noether-Fano
inequality we may assume that $D$ is a prime divisor. According
to Proposition 3.1 which we have just proven, $\nu\leq\frac32n$,
so that we are in the situation of Sec. 2.2. The pair
$(V^+,\frac{1}{n}D^+)$ is not log canonical, some subvariety
$S\subset E$ is the centre of a non log canonical singularity of
that pair. We assume that $S$ has the maximal dimension among all
centres of such singularities, so that $\mathop{\rm dim} S\leq 4$.

{\bf Proposition 3.2.} {\it $S$ is of codimension at least two in
the exceptional quadric $E$.}

{\bf Proof.} If $S\subset E$ is a prime divisor, then, in
accordance with Proposition 2.2, $S$ is a hyperplane section of
the quadric $E$. Let $P\ni o$ be the unique hyperplane in
${\mathbb P}$, cutting out $S$ on $E$, that is,
$$
V^+_P\cap E=S,
$$
where $V_P=V\cap P$, $V^+_P\subset V^+$ is the strict transform.
The pair
$$
(V,V_P)
$$
is canonical, so that $D\neq V_P$ and the set-theoretic
intersection $D\cap V_P$ is of codimension two. For the effective
cycle $D_P=(D\circ V_P)$ of codimension two we get
$$
\mathop{\rm mult}\nolimits_oD_P=\mathop{\rm
mult}\nolimits_oD+2\mathop{\rm mult}\nolimits_SD^+>4n,
$$
so that, arguing as in the proof of Proposition 3.1, we construct
a sequence of irreducible subvarieties
$$
Y_2,Y_3,\dots,Y_{M-1},
$$
$\mathop{\rm codim}Y_i=i$, $Y_2$ is an irreducible component of
the cycle $D_P$, satisfying the inequality
$$
\frac{\mathop{\rm mult}_o}{\mathop{\rm
deg}}Y_2\geq\frac{\mathop{\rm mult}_o}{\mathop{\rm
deg}}D_P>\frac{4}{M+1},
$$
and obtain a contradiction:
$$
\frac{\mathop{\rm mult}_o}{\mathop{\rm
deg}}Y_{M-1}>\frac{4}{M+1}\cdot\frac54\cdot\dots\cdot\frac{M+1}{M}=1,
$$
which is impossible. Q.E.D. for the proposition.

We conclude that the second case of Proposition 2.2 takes place:
there is a hyperplane section $Z\supset S$ of the exceptional
quadric $E$, satisfying the inequality (\ref{b3}). Let
$P\subset{\mathbb P}$ be the unique hyperplane, cutting out $Z$ on
$E$ (in the same sense as in the proof of Proposition 3.2),
$V_P=V\cap P\neq D$. For the effective cycle $D_P=(D\circ V_P)$
of codimension two we get
$$
\mathop{\rm mult}\nolimits_oD_P\geq\mathop{\rm
mult}\nolimits_oD+2\mathop{\rm
mult}\nolimits_ZD^+>\frac43(l+n)>\frac83n.
$$
Unfortunately, this estimate is insufficient for excluding the
maximal singularity in the same way which we used in the proof of
Propositions 3.1 and 3.2. (The lower bound for $(\mathop{\rm
mult}_o/\mathop{\rm deg}Y_{M-1}$ turns out to be smaller than
one, which does not allow to get a contradiction.) However, we
can consider the pair
$$
(V_P,\frac{1}{n}D_P).
$$
By inversion of adjunction its strict transform
$$
(V^+_P,\frac{1}{n}D^+_P)
$$
with respect to the blow up of the point $o$ is not log
canonical, whereas the subvariety $S\subset E_P=Z$ is a centre of
a non log canonical singularity of that pair. We may assume that
$$
\mathop{\rm mult}\nolimits_oD^+_P\leq 4n,
$$
otherwise we could obtain a contradiction using word for word the
same argument as in the proof of Proposition 3.2. Without loss of
generality we may assume that $S$ is a maximal centre of a non
log canonical singularity of the pair $(V^+_P,\frac{1}{n}D^+_P)$.
By Remark 2.2 (and the inequality $\mathop{\rm dim}S\geq 4$) we
may apply Proposition 2.2 to the latter pair and obtain the
alternative: either

1) $S$ is a hyperplane section of the quadric $E_P$, or

2) there is a hyperplane section $Z^+\supset S$ of the quadric
$E_P$, satisfying the inequality
\begin{equation} \label{c1}
\mathop{\rm mult}\nolimits_{Z^*}D^+_P>\frac{2n-l^*}{3},
\end{equation}
where $D^+_P\sim nH_P-l^*E_P$. Recall that the integer $l^*$
satisfies the inequality
$$
l^*>\frac23(l+n)>\frac43n.
$$
\vspace{0.3cm}

%%%%%%%%%%%%%%%%%%%%%%%%%%%%%%%%%%%%%%%%%%%%%%%%%%%%%%%%%%%%%%%%%%%%%
%%%%%%%%%%%%%%%%%%%%%%%%%%%%%%%%%%%%     subsection 3.3

{\bf 3.3. The repeated hyperplane section.} Let $R\subset
P={\mathbb P}^M$ be the unique hyperplane, cutting out on $E_P$
in the case 1) the subvariety $S$, in the case 2) the subvariety
$Z^*$.

Assume that the case 1) takes place. By linearity of the
inequalities
\begin{equation}\label{c2}
\mathop{\rm
mult}\nolimits_oD_P>\frac83n\quad\mbox{and}\quad\mathop{\rm
mult}\nolimits_SD^+_P>n
\end{equation}
in the divisor $D_P$ and by the fact that the divisor
$V_R=V_P\cap R$ does not satisfy them, we may assume that the
divisor $D_P$ (which is, possibly, reducible) does not contain
$V_R$ as a component (otherwise, we delete this component, which
could only make both inequalities (\ref{c2}) sharper). For this
reason, the intersection $D_P\cap V_R$ is of codimension two on
$V_P$ and the effective algebraic cycle
$$
D_R=(D_P\circ V_R)
$$
is well defined and satisfies the inequality
\begin{equation}\label{c3}
\mathop{\rm mult}\nolimits_oD_R>\mathop{\rm
mult}\nolimits_oD_P+2\mathop{\rm
mult}\nolimits_SD^+_P>\frac{14}{3}n.
\end{equation}
By linearity of the last inequality we may assume that $D_R=Y$ is
an irreducible variety, that is, a prime divisor on $V_R$.
However, $V_R$ is a section of the hypersurface $V$ by the linear
subspace $R\subset{\mathbb P}$ of codimension two. Let
$D_2\,|\,_R=\{q_2\,|\,_{V_R}=0\}$ be the first hypertangent
divisor of the variety $V_R$. According to the regularity
condition (R2.2), $D_2\,|\,_R$ is irreducible and does not
satisfy the inequality (\ref{c3}), that is,
$$
Y\neq D_2\,|\,_R
$$
and for that reason $Y\not\subset D_2=\{q_2\,|\,_V=0\}$. Let
$Y_4$ be an irreducible component of the effective algebraic cycle
$(Y\circ D_2)$ with the maximal ratio $\mathop{\rm
mult}_o/\mathop{\rm deg}$. We have the inequality
$$
\frac{\mathop{\rm mult}\nolimits_o}{\mathop{\rm
deg}}Y_4>\frac{7}{M+1}.
$$
Now we argue as above: construct a sequence of irreducible
subvarieties
$$
Y_4,Y_5,\dots,Y_{M-1},
$$
$\mathop{\rm codim}_VY_i=i$, $Y_{i+1}$ is an irreducible
component of the effective cycle $(Y_i\circ D_{i+2})$, where
$D_{i+2}\in\Lambda_{i+2}$ is a generic hypertangent divisor. For
$Y_{M-1}$ we get the estimate
$$
\frac{\mathop{\rm mult}\nolimits_o}{\mathop{\rm
deg}Y_{M+1}}>\frac{7}{M+1}{\cdot\frac{M+1}{6}}=\frac76,
$$
which is impossible. The contradiction excludes the case 1).

Finally, consider the hardest case 2). Again we use the linearity
of the conditions that are satisfied by the divisor $D_P$ on
$V_P$: the inequality
$$
\mathop{\rm mult}\nolimits_oD_P>\frac83n
$$
and the existence of a non log canonical singularity of the pair
$(V^+_P,\frac{1}{n}D^+_P)$ with the centre $S$. As above, set
$V_R=V_P\cap R$. Since $\mathop{\rm mult}_oV_R=2<\frac83$ and the
pair $(V^+_P,V_R)$ is log canonical, we may assume that the
divisor $D_P$ does not contain the hyperplane section $V_R$ as a
component. For this reason, the effective cycle of codimension two
(with respect to $V_P$) $D_R=(D_P\circ V_R)$ is well defined.
This cycle satisfies the inequality
\begin{equation}\label{c4}
\mathop{\rm mult}\nolimits_oD_R\geq\mathop{\rm
mult}\nolimits_oD_P+\mathop{\rm
mult}\nolimits_{Z^*}D^+_P>\frac{28}{9}n.
\end{equation}
As in the case 1), the inequality (\ref{c4}) is sufficient to
conclude (by the regularity condition (R2.2)), that the component
$Y=Y_3$ of the cycle $D_R$ with the maximal ratio $\mathop{\rm
mult}_o/\mathop{\rm deg}$ is not contained in the divisor
$D_2=\{q_2\,|\,_V=0\}$, so that the effective cycle
$$
(Y_3\circ D_2)
$$
is well defined. There is an irreducible component $Y_4$ of this
cycle, satisfying the inequality
\begin{equation}\label{c5}
\frac{\mathop{\rm mult}\nolimits_o}{\mathop{\rm
deg}}Y_4>\frac{14}{3(M+1)}.
\end{equation}
By the regularity conditions (R2.2) this procedure can be
repeated three times more. Let us consider in more details the
first step. The subvariety
$$
W_{2\cdot 3}=\{q_2=0\}\cap\{q_3=0\}\cap V_R
$$
is irreducible, of degree $\mathop{\rm deg}W_{2\cdot 3} =6(M+1)$
and multiplicity
$$
\mathop{\rm mult}\nolimits_oW_{2\cdot 3} =24,
$$
so that $W_{2\cdot 3} \neq Y_4$. Since by construction
$$
Y_4\subset\{q_2=0\}\cap V_R,
$$
this implies that $Y_4\not\subset\{q_3=0\}$, so that
$$
Y_4\not\subset D_3=\{(q_2+q_3)\,|\,V=0\}
$$
and the effective cycle $(Y_4\circ D_3)$ of codimension 5 is well
defined. Some irreducible component $Y_5$ of this cycle satisfies
the inequality
$$
\frac{\mathop{\rm mult}\nolimits_o}{\mathop{\rm
deg}}Y_5>\frac{56}{9(M+1)}.
$$
In the same way, consider the irreducible subvariety
$$
W_{2\cdot 3\cdot 4}=W_{2\cdot 3}\cap\{q_4=0\}
$$
and construct an irreducible subvariety $Y_6\subset V$ of
codimension 6, satisfying the inequality
$$
\frac{\mathop{\rm mult}\nolimits_o}{\mathop{\rm
deg}}Y_6>\frac{70}{9(M+1)}.
$$
Finally, consider the subvariety
$$
W_{2\cdot 3\cdot 4\cdot 5}=W_{2\cdot 3\cdot 4}\cap\{q_5=0\}
$$
and construct an irreducible subvariety $Y_7\subset V$ of
codimension 7, satisfying the inequality
$$
\frac{\mathop{\rm mult}\nolimits_o}{\mathop{\rm
deg}}Y_7>\frac{28}{3(M+1)}.
$$
If $M=8$, we get a contradiction. If $M\geq 9$, we apply the
technique of hypertangent divisors in the traditional version (as
in the proof of Proposition 3.1), intersecting $Y_7$ with the
generic hypertangent divisors $D_9,\dots$. For the irreducible
curve $Y_{M-1}$ we obtain the estimate
$$
\frac{\mathop{\rm mult}\nolimits_o}{\mathop{\rm
deg}}Y_{M-1}>\frac{28}{3(M+1)}\cdot\frac{M+1}{9}=\frac{28}{27}.
$$
This contradiction completes the case 2) which is now excluded.
Q.E.D. for Theorem 5.

%%%%%%%%%%%%%%%%%%%%%%%%%%%%%%%%%%%%%%%%%%%%%%%%%%%%%%%%%%%%%%%%%
%%%%%%%%%%%%%%%%%%%%%%%%%%%%%%%%%%%%%%%%%%%%%%%%%%%%%%%%%%%%%%%%%
%%%%%%%%%%%%%%%%%%%%%%%%%%%%%%%%%%% references

{\small

\section*{References}

\noindent 1. Pukhlikov A.V., Birational geometry of Fano direct
products, Izvestiya: Mathematics, V. 69 (2005), no. 6, 1225-1255.
\vspace{0.3cm}

\noindent 2. Pukhlikov A.V., Birationally rigid varieties. I.
Fano varieties. Russian Math. Surveys. 2007. V. 62, no. 5,
857-942. \vspace{0.3cm}

\noindent 3. Pukhlikov A.V., Birationally rigid Fano
hypersurfaces with isolated singularities, Sbornik: Mathematics
{\bf 193} (2002), No. 3, 445-471. \vspace{0.3cm}

\noindent 4. Pukhlikov A.V., Birational automorphisms of Fano
hypersurfaces, Invent. Math. {\bf 134} (1998), no. 2, 401-426.
\vspace{0.3cm}

\noindent 5. Cheltsov I.A., Fano varieties with many self-maps.
Adv. Math. {\bf 217} (2008), no. 1, 97-124. \vspace{0.3cm}

\noindent 6. Koll{\'a}r J., et al., Flips and Abundance for
Algebraic Threefolds, Asterisque 211, 1993. \vspace{0.3cm}

\noindent 7. Cheltsov I.A., Log canonical thresholds on
hypersurfaces. Sbornik: Mathematics. 2001. V. 192, no. 7-8,
1241-1257. \vspace{0.3cm}}

\end{document}